\documentclass[11pt, reqno]{amsart}
\usepackage{amsmath, amsthm, amsfonts, amssymb, enumerate, dsfont, graphicx, epstopdf}
\usepackage{float}
\usepackage{hyperref,natbib}
\hypersetup{
  colorlinks=true,        % false: boxed links; true: colored links
  citecolor=green,        % color of links to bibliography
}
\usepackage[hmargin=3cm, vmargin=3cm]{geometry}

\newtheorem{thm}{Theorem}[section]
\newtheorem{cor}[thm]{Corollary}

\theoremstyle{definition}
\newtheorem{defn}[thm]{Definition}
\newtheorem{prob}{Problem}
\newtheorem{notn}{Notation}[section]
\newtheorem{eg}{Example}[section]
\theoremstyle{remark}
\newtheorem{qstn}{Question}
\newtheorem{remark}{Remark}

\numberwithin{equation}{section}

\def\beq{\begin{equation}}
\def\eeq{\end{equation}}
\def\ba{\begin{enumerate}[(a)]}
\def\bei{\begin{enumerate}[(i)]}
\def\be{\begin{enumerate}[(1)]}
\def\ee{\end{enumerate}}
\def\bi{\begin{itemize}}
\def\ei{\end{itemize}}
\def\bd{\begin{defn}}
\def\ed{\end{defn}}
\def\bqn{\begin{qstn}}
\def\eqn{\end{qstn}}
\def\bt{\begin{thm}}
\def\et{\end{thm}}
\def\bc{\begin{cor}}
\def\ec{\end{cor}}
\def\bark{\begin{remark}}
\def\eark{\end{remark}}
\def\bp{\begin{prop}}
\def\ep{\end{prop}}
\def\bpb{\begin{prob}}
\def\epb{\end{prob}}
\def\beg{\begin{eg}}
\def\eeg{\end{eg}}
\def\bark{\begin{remark}}
\def\eark{\end{remark}}
\def\bn{\begin{notn}}
\def\en{\end{notn}}
\def\bpf{\begin{proof}}
\def\epf{\end{proof}}
\def\RR{\ensuremath{\mathbb{R}}}

\def\LP{\left\{} 
\def\RP{\right\}}

\DeclareMathOperator{\Exp}{\mathbb{E}}

\usepackage{bigstrut}
\usepackage{booktabs}

\begin{document}

\title[Sparse Minimax Optimality of Bayes PDE from Clustered Discrete Priors]{Sparse Minimax Optimality of Bayes Predictive Density Estimates from Clustered Discrete Priors}
\date{}
\author{Ujan Gangopadhyay}
\address{Ujan Gangopadhyay, \ Department of Mathematics, \ University of Southern California, \ Los Angeles, CA, USA.}
\email{ujan.gangopadhyay@usc.edu}
\author{Gourab Mukherjee}
\address{Gourab Mukherjee, \ Marshall School of Business, \ University of Southern California, \ Los Angeles, CA, USA.}
\email{gourab@usc.edu}

\begin{abstract}
We consider the problem of predictive density estimation under Kullback-Leibler loss in a high-dimensional Gaussian model with exact sparsity constraints on the location parameters. We study the first order asymptotic minimax risk of Bayes predictive density estimates based on product discrete priors where the proportion of non-zero coordinates converges to zero as dimension increases. Discrete priors that are product of clustered univariate priors provide a tractable configuration for diversification of the future risk and are used for constructing efficient predictive density estimates. We establish that the Bayes predictive density estimate from an appropriately designed clustered discrete prior is asymptotically minimax optimal. The marginals of our proposed prior have infinite clusters of identical sizes. The within cluster support points are equi-probable and the clusters are periodically spaced with geometrically decaying probabilities as they move away from the origin. The cluster periodicity depends on the decay rate of the cluster probabilities. Under different sparsity regimes, through numerical experiments, we compare the maximal risk of the Bayes predictive density estimates from the clustered prior with varied competing estimators including those based on geometrically decaying non-clustered priors of \citet{Johnstone94a} and \citet{mukherjee2017minimax} and obtain encouraging results.  
\end{abstract}



\maketitle

\bibliographystyle{plainnat}

\section{Introduction and Main Results}
A fundamental problem in statistical prediction analysis is to choose a probability distribution based on observed data that will be good in predicting the behavior of future samples \citep{Aitchison-book,Geisser-book,Aitchison75}. The future probability density conditioned on the observed past is referred to as the predictive density and estimating it plays an important role in a number of statistical applications \citep{Liang-thesis,Mukherjee-thesis}. 
Consider the problem of predictive density estimation in a $n$-dimensional Gaussian location model where the observed past vector 
$X\sim N_n(\theta,v_xI)$ and the future vector $Y\sim N_n(\theta,v_yI)$. The variances $v_x$ and $v_y$ are known. The future and past vectors are related only through the unknown location vector $\theta$. Consider predictive density estimators (\textit{prde}) $\hat{p}(y|x)$ and measure their performance in estimating the true future density $p(y|\theta,v_y)=N_n(\theta,v_yI)$ by the global divergence measure of \citet{Kullback51}, 
\begin{align}\label{eq:loss}
L(\theta,\hat{p}(\cdot|x))=\int p(y|\theta,v_y)\log\bigg(\frac{p(y|\theta,v_y)}{\hat{p}(y|x)}\bigg)\,dy.
\end{align}
The KL risk integrates the above KL loss over the past distribution and is given by
%\[\label{eq:risk}
$\rho(\theta,\hat{p})=\int L(\theta,\hat{p}(\cdot|x))p(x|\theta,v_x)\,dx.$
%\]
Given any prior $\pi$ on $\theta$, the Bayes \textit{prde} 
$\hat{p}_\pi(y|x)=\int p(y|\theta,v_y)\pi(d\theta|x)$. The average integrated risk
$B(\pi,\hat{p})=\int\rho(\theta,\hat{p})\pi(d\theta)$, when well-defined, is minimized by $\hat{p}_\pi$ yielding the 
Bayes risk $B(\pi)=\inf_{\hat{p}}B(\pi,\hat{p})$. 

As dimension $n$ increases, there exists decision theoretic parallels between \textit{prde} under \eqref{eq:loss} and point estimation (\text{PE}) of the multivariate normal mean under square error loss (see \citealp{George06,George12,Komaki01,Fourdrinier11,maruyama2016harmonic,kubokawa2013minimaxity,ghosh2018hierarchical,Xu10,Brown08,Ghosh08}). Sparse \textit{prde} under exact $\ell_0$ sparsity constraints on the location parameter is studied in \citet{mukherjee2017minimax,Mukherjee-15} where efficacy of different \textit{prde}s were evaluated with respect to the minimax benchmark risk 
$R^*(\Theta)=\inf_{\hat{p}}\sup_{\theta\in\Theta}\rho(\theta,\hat{p})$.
For an $\ell_0$ constrained parameter space
$\Theta_0[s_n]=\LP\theta\in\RR^n: \sum_{i=1}^n 1\{\theta_i \neq 0 \}\leq s_n\RP$ when
$\eta_n=s_n/n \to 0$, the first order asymptotic minimax risk was evaluated as
\[
R^*(\Theta_0[s_n])=(1+r)^{-1} n\, \eta_n\log\eta_n^{-1} (1+o(1)) \text{ as } n \to \infty,
\]
where $r=v_y/v_x$. The minimax risk increases as $r$ decreases. The difficulty of the density estimation problem increases as $r$ decreases as we need to estimate the future observation density based on increasingly noisy past observations. The rate of convergence of the minimax risk with $n$ does not depend on $r$, and so exact determination of the constants is needed to show the role of $r$ in this prediction problem. Several predictive phenomena that contrast with point estimation results have been reported with the divergence becoming palpable as $r$ decreases. 

Here, we study the risk of Bayes predictive density estimators based on sparse discrete priors. In order to incorporate the knowledge on sparsity of the parameters,  we consider priors with an atom of probability (spike) at the origin. Spike-and-slab priors based procedures have been shown to be very successful for sparse estimation \citep{Johnstone04,clyde2000flexible,rovckova2018spike}. Here, we consider slabs based on periodic discrete priors. 
%In several minimax problems, particularly for constrained inferences, the least favorable prior is often discrete \citep{bickel1983minimax,Johnstone-book,marchand2004estimation} and 
Risk analysis of estimators based on discrete priors has a rich history in statistical decision theory \citep{Johnstone-book,marchand2004estimation}, particularly for studying the worst-case geometry of parametric spaces \citep{bickel1983minimax,kempthorne1987numerical}. 
\citet{Johnstone94a} (henceforth referred to as J94) established that for sparse point estimation a product prior based on discrete marginals containing equi-spaced support-points with geometrically decaying probability  is asymptotically minimax optimal. \citet{mukherjee2017minimax} (referred hereon as MJ17) showed that Bayes \textit{prdes} from such grid priors are minimax sub-optimal. The clustered discrete prior we study here is inspired by the risk diversification phenomenon introduced in \citet{Mukherjee-15} (referred to as MJ15) for constructing minimax optimal \textit{prde}s. MJ15 showed that in contrast to point estimation, for obtaining minimax optimality in sparse \textit{prde} we need to incorporate the notion of diversification of the future risk. A product prior consisting of clustered discrete marginals with equi-probable support points in each clusters were used along with thresholding. Here, we conduct detailed worst-case risk analysis of \textit{prde}s based on generic versions of such clustered discrete priors. As such, MJ15 used a version of the Bayes \textit{prdes} that was based on only the origin adjoining  two clusters of the prior analyzed here. Our proposed clustered prior based Bayes \textit{prde} also has the advantage of avoiding the discontinuous thresholding operation in order to obtain sparse minimax optimality. The risk analysis of predictors based on clustered priors differs in fundamental aspects from the analysis of non-clustered priors in MJ17 and provides new insights on the risk profiles of segmented priors. %which might be of independent interest. 
Next, we present our main result following which detailed background and connections to the existing literature is provided. 
\begin{table}[t]
	\begin{center}
		\caption{The size $K_r$ of each cluster in our proposed univariate cluster prior $\pi_{\sf C}$ as $r$ varies.}
		\scalebox{0.7}{	\begin{tabular}{|c|c|c|c|c|c|c|c|c|}
				\hline
				$r$     & 0.0654 & 0.0759 & 0.0910 & 0.1150 & 0.1601 & 0.2826 & 0.5000 & \multicolumn{1}{l|}{$>$0.5000} \\
				\hline
				$K_r$     & 8     & 7     & 6     & 5     & 4     & 3     & 2     & 1 \bigstrut\\
				\hline
			\end{tabular}%
		}
		\label{ref:tab-sf}%
	\end{center}
\end{table}%

\begin{figure}[h]
	% The arguments in the next line are {height}{optional width}{used only by OUP for typesetting}[filename, in directory art]
	\includegraphics[width=\textwidth]{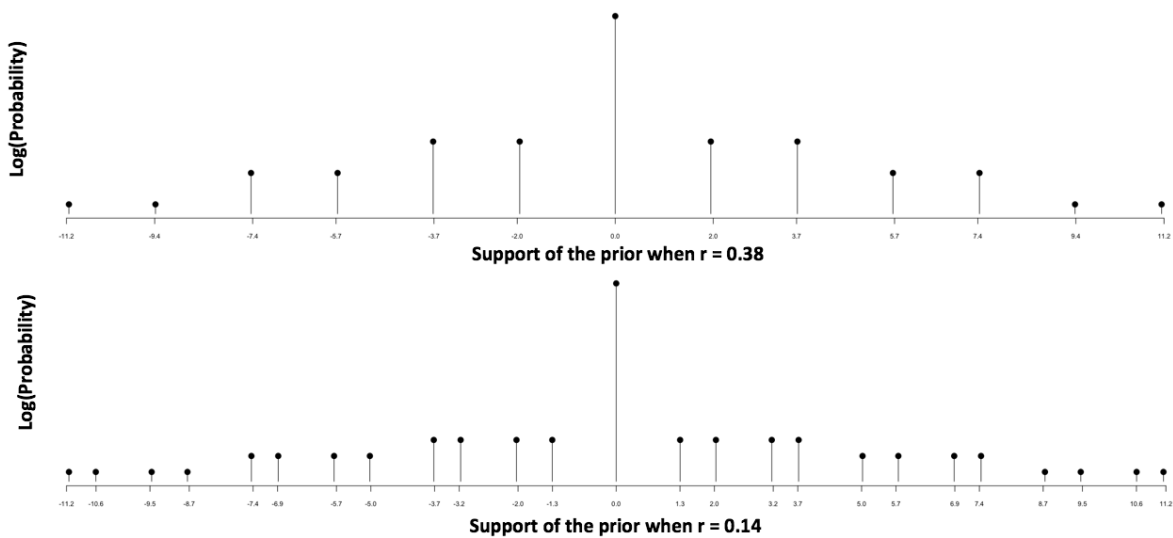}
	% note that files may not be rotated
	\caption{Schematic for our proposed univariate cluster prior when $r$ equals 0.38 (top) and 0.14 (bottom) respectively. The x-axis shows the spacings between and within the clusters and the y-axis the logarithm of the prior probabilities. Figure drawn to scale with $\eta=0.001$. Only the six clusters are displayed with the rest being truncated. }
	\label{fig1}
\end{figure}

\textit{Main Result.}
%Here, we study asymptotic minimaxity of Bayes \textit{prdes  based  on discrete product priors with symmetric marginals. 
%Consider a generic class of univariate symmetric priors:
%$$\pi_{\sf C}[\eta,r]=(1-\eta)\delta_0+\frac{1-\eta}{2} \sum_{i=1}^{\infty} \eta^i \big\{C_i (\eta,r) + C_{-i}(\eta,r)\big\}$$
For any fixed positive $r$, consider the Bayes \textit{prde} from a discrete product prior consisting of symmetric marginals $\pi_{\sf CL}$ (defined below). The marginal has equi-spaced clusters of atoms with geometrically decaying probability content in the clusters as they move away from the origin. For any $\eta \in (0,1)$ and $r \in (0,\infty)$ consider the univariate clustered discrete prior: 
\begin{align}\label{eq:cl-1}
\pi_{\sf CL}[\eta,r; \gamma,\kappa]=(1-\eta)\delta_0+\frac{1-\eta}{2} \sum_{i=1}^{\infty} \eta^{i} \big\{C_i (\eta,r; \gamma, \kappa) + C_{-i}(\eta,r; \gamma, \kappa)\big\}~,
\end{align}
which has an atom of probability $1-\eta$ at the origin and the remaining $\eta$ probability shared across clusters. 
Each of the clusters $C_i$ has $\kappa$ atoms $\{\mu_{ij}: j=1,\ldots,\kappa\}$ of equal probability which is the reason for referring such prior distributions as \textit{clustered priors}. Let  $v=(1+r^{-1})^{-1}$, $\lambda_e:=\lambda_e(\eta)=(-2v_x\log\eta)^{1/2}$ and $\lambda_f:=\lambda_f(\eta,r)=v^{1/2}\lambda_e$.
For any fixed $\gamma \geq 1$, the atoms in $C_1$ are aligned in between $\lambda_f$ and $\lambda_e$ in a geometric progression with common ratio $\gamma$, i.e.,   
$\mu_{1j}(\eta,r,\gamma)= \gamma^{j-1} \lambda_f  \wedge \lambda_e$ for $1\leq j \leq \kappa$. 
Such geometric spacing was introduced in MJ15 (see Theorem 1C)
For $i\geq 2$ the atoms are extended  periodically to cluster $C_{i}$ as $\mu_{ij}= (i-1) \mu_{1\kappa} + \mu_{1j}$ and by symmetry $ \mu_{-ij}=-\mu_{ij}$ to the negative axis.  
Thus, the clusters themselves are equidistant at a separation of $\lambda_f$ and while the atoms within each cluster has equal probability, the clusters themselves have geometrically decaying probabilities: 
\begin{align}\label{eq:cl-2}
C_i(\eta,r; \gamma, \kappa)= \frac{1}{\kappa} \sum_{j=1}^{\kappa}  \delta_{\mu_{ij}}  \text{ and } P(C_i)= 2^{-1}(1-\eta)\eta^{|i|} \text{ for } i \in \mathbb{Z}\setminus \{0\}.
\end{align}
Our proposed cluster prior $\pi_{\sf C}$ has $\gamma=\gamma_r$ and $\kappa=K$ where, $\gamma_r=1+4 r$ and
\begin{align}\label{eq.K}
K:=K_r=1 +  \big \lceil \log(1+r^{-1})/(2\log \gamma_r)\big \rceil  \cdot 1\{r< r_0\}~.
\end{align}
Thus, 
$\pi_{\sf C}[\eta,r]:=\pi_{\sf CL}[\eta,r;\gamma_r,K]$. Here, $r_0 = 0.5$. Note that, $K=1$ iff $r \geq r_0$. The significance of $r_0$ is shown in Proposition 1 of the supplementary materials. When $K \geq 3$ and $i \geq 1$, all atoms except the $K$th one in any cluster $C_i$ are aligned in a geometric progression starting from $\mu_{i1}= (i-1) \lambda_e + \lambda_f$, with common ratio $1+4r$ and $\mu_{iK}=i\lambda_e$.  Table~\ref{ref:tab-sf} shows the cluster size as $r$ varies. Figure~\ref{fig1} shows the schematic diagram of the (truncated) prior with 6 clusters for two instances when $r=0.38$ and $r=0.14$ respectively. While the former has clusters of size 2, the latter has cluster size 4. Figure~\ref{fig1} illustrates a key aspect of the cluster prior: for $r< r_0$ the gap $\mu_{i,K}-\mu_{i,K-1}$ is allowed to vary widely with $r$ while $\mu_{i+1,1}-\mu_{i,K}$ is fixed at $\lambda_f$ for all $i$.  
%The common ratio here is double of that of MJ15 and hence the number of cluster points (compared to table 1 of MJ15) has been reduced for smaller $r$ values.
% Table generated by Excel2LaTeX from sheet 'Sheet3'
%$\pi_n(\theta)=\prod_{i=1}^n \pi[]$
% on independent univariate  mixture priors $\pi_{\sf C}[\eta]$ with an atom of probability at $0$ and a symmetric discrete prior consisting of periodic clusters. Consider the univariate prior: 

Now, consider the multivariate clustered prior $\pi^{\sf C}_n[\eta_n,r](d\theta)= \prod_{i=1}^n \pi_{\sf C}[\eta_n,r](d\theta_i)$ on $\mathbb{R}^n$. %$\Theta_0[s_n]$. 
Then, the Bayes \textit{prde} $\hat p_{\sf C}[\eta_n,r]$ based on $\pi^{\sf C}_n[\eta,n]$ is asymptotically minimax optimal.

%\[
%r=v_y/v_x=v_y\quad\mbox{and}\quad .
%\]

\bt Fix any $r\in(0,\infty)$. If $\eta_n=s_n/n \to 0$, then
\[
\lim_{n \to \infty} \; \bigg\{\sup_{\theta \in \Theta_0[s_n]} \rho\big(\theta,\hat p_{\sf C}[\eta_n, r]\big)\bigg\}\bigg/ R^*(\Theta_0[s_n])=1.
\]
%Additionally, if $s_n \to \infty$, then $\pi^{\sf C}$ is also asymptotically least favorable: 
%$$\liminf_{n \to \infty}  B(\pi^{\sf C}_n)/ R^*(\Theta_0[s_n])=1.$$
\et

\textit{Background.} For understanding the decision theoretic implications of the above result, we briefly revisit the risk properties of sparse product priors based on symmetric marginals.  
%of the form 
%$(1-\eta)\delta_0+\frac{\eta}{2} (\gamma_B^+ + \gamma_B^-) $
It follows from J94  that for point estimation of the normal mean over $\Theta_0[s_n]$ under $\ell_2$ loss, the posterior mean of the grid prior $\pi_n^{\sf EG}$ is minimax optimal as $\eta_n \to 0$. $\pi_n^{\sf EG}$ constitutes of i.i.d. copies of univariate grid prior $\pi_{\sf EG}[\eta_n,r]$ which is defined for any fixed $r$ and $\eta \in (0,1)$ as
$$\pi_{\sf EG}[\eta,r]=(1-\eta)\delta_0+\frac{1-\eta}{2} \sum_{i=1}^{\infty} \eta^i \big\{\delta_{i\lambda_e} + \delta_{-i\lambda_e}\big\}~.$$
In contrast to $\pi_{\sf C}$, $\pi_{\sf EG}$ always has only one point in each cluster. However, they have identical probability decay rate as the clusters extend away from the origin. 
MJ17 showed that the \textit{prde} based on $\pi_n^{\sf EG}$ is sub-optimal for \textit{prde} estimation based on KL loss. The Bayes \textit{prde} based on a product grid prior whose univariate marginals $\pi_{\sf PG}$ 
(subscripts PG and EG denote predictive and estimative grids)  has reduced spacing between the atoms and reduced probability decay rate, was established to be minimax optimal in the predictive regime abet for $r \geq \tilde r_0=(\sqrt 5 - 1)/4$:
$$\pi_{\sf PG}[\eta,r]=(1-\eta)\delta_0+\frac{\eta(1-\eta^v)}{2} \sum_{i=1}^{\infty} \eta^{(i-1)v} \big\{\delta_{i\lambda_f} + \delta_{-i\lambda_f}\big\}~.$$
For constructing a minimax optimal Bayes \textit{prde} for all values of $r$, MJ17 suggested using a bi-grid prior with two different sections: inner and outer. While the outer section has the spacing and decay rate of $\pi_{\sf PG}$ the inner section has further reduced spacing.  Let $b:=b(r)=\min\{4r(1+r)/(1+2r),1\}$ and $J=1+\lceil 2 b^{-3/2}\rceil$. For any  integer $j$ and $l$, define the inner section support points ${\sf I}_j=\text{sign}(j)\{\lambda_f+b(|j|-1)\lambda_f\}$ 
and the outer section atoms ${\sf O}_l=\text{sign}(l)\{ I_J+ |l| \lambda_f\}$. Then, the univariate bi-grid prior is: 
$$\pi_{\sf BG}[\eta,r]=(1-\eta)\delta_0+\frac{\eta \, c(\eta,r)}{2} \bigg[\sum_{j=1}^{J} \eta^{(j-1)b^2v} \big\{\delta_{{\sf I}_j} + \delta_{{\sf I}_{-j}}\big\} + \eta^{(J-1)b^2v}
 \sum_{l=1}^{\infty} \eta^{lv} \big\{\delta_{{\sf O}_l} + \delta_{{\sf O}_{-l}}\big\} \bigg]
$$
where, $c(\eta,r)$ is the normalizing constant defined in eqn. (28) of MJ17. The multivariate prior $\prod_{i=1}^n \pi_{\sf BG}[\eta_n,r] (d\theta_i)$ is minimax optimal for any $r$. 
Note that $\pi_{\sf BG}$ agrees with $\pi_{\sf PG}$ for $r \geq \tilde r_0$. 

\textit{Discussion.} Unlike the univariate grid priors $\pi_{\sf EG}, \pi_{\sf PG}$ where support points has geometric probability decay, $\pi_{\sf C}$ has support points with identical probability within each clusters. The clusters in $\pi_{\sf C}$ however has the same decay rate as the support points in $\pi_{\sf EG}$. The maximum gap between atoms in $\pi_{\sf C}$ equals the spacing in $\pi_{\sf PG}$. Equiprobable atoms in the clusters was introduced in MJ15 to control predictive risk via the new notion of risk diversification. As such consider a truncated cluster prior with only two clusters:
$\pi_{\sf TC}[\eta,r]=(1-\eta)\delta_0+{\eta}/{2}  \{C_1  + C_{-1}\}$ where $C_1=C_1(\eta,r;\tilde \gamma_r ,\tilde K_r)$ as in \eqref{eq:cl-2}
with $\tilde \gamma_r=1+2r$ and $\tilde K_r$ given by $K_r-1$ with the formula in \eqref{eq.K} used with $\tilde \gamma_r$ in place of $\gamma_r$. As the prior $\pi_{\sf TC}$ is bounded at $\lambda_e$, its corresponding Bayes \textit{prde} $\hat p_{\sf CT}$ has unbounded risk. Thresholded product \textit{prde} $\hat p^{\sf T}_n(y|x)=\prod_{i=1}^n\hat p_{\sf T}(y_i|x_i)$ with
$$\hat p_{\sf T} (y_i|x_i) = \hat p_{\sf TC}[\eta_n, r](y_i|x_i) 1\{|x_i|\leq \lambda_e(\eta_n)\} + \phi(y_i|x_i, v_x+v_y) 1\{|x_i|> \lambda_e(\eta_n)\}~$$
was shown in MJ15 to be minimax optimal. 
Note that, the thresholding was done at the boundary $\lambda_e(\eta_n)$ of the truncated univariate prior; above the threshold the Bayes \text{prde} based on the uniform prior, which is Gaussian with variance $v_x+v_y$, was used. Thresholding rules are not smooth functions of the data and it was conjectured in Sec. 6 of MJ15  that periodic clustered priors of the form of \eqref{eq:cl-1}-\eqref{eq:cl-2} can attain minimax optimality without the discontinuous thresholding operation. Here, we study the risk properties of such cluster priors and establish minimax optimality of the properly calibrated prior $\pi_C$. We found that the common ratio $\tilde \gamma_r$ used in MJ15 was not optimal and can be increased to $\gamma_r$. However, as a consequence of removing thresholding  we needed one more atom than MJ15 in our proposed cluster prior $\pi_{\sf C}$ for small values of $r$. 

The new phenomenon of risk diversification introduced in MJ15 to obtain minimax optimality of \textit{prde}s was further extended in MJ17 where it was shown that to attain minimax optimality of Bayes \textit{prde}s based on discrete priors, the atoms need to be much denser near the origin that away from the origin. The inner section spacing $b(r)$ of the  bi-grid prior $\pi_{\sf BG}$ of MJ17 is slightly lower but quite close to the minimal within cluster spacing in $\pi_{\sf C}$.  An intrinsic difference between $\pi_{\sf C}$ and $\pi_{\sf BG}$ is that for $\eta \to 0$ the first cluster $C_1$ protrudes much beyond  inner section of $\pi_{\sf BG}$, particularly for smaller values of $r$.  
Though the Bayes \textit{prde}s from the cluster prior and the bi-grid prior are both minimax optimal (compare theorem 1  here with theorem 1.2 of MJ17), there exists interesting disparity in geometry of their manifolds; subsequently, their maximal risk for them are controlled by different facets of the risk diversification principle. This necessitates separate analysis and proofs of the risk properties of $\pi_{\sf C}$ than that of bi-grid priors.  %In lieu of these differences,   demonstrate that there exists interesting disparity in geometry of the manifolds of different multivariate sparse discrete priors whose Bayes \textit{prde} are asymptotically minimax optimal over $\Theta_0(s_n)$ as $\eta_n \to 0$. 
\begin{figure}[h]
	% The arguments in the next line are {height}{optional width}{used only by OUP for typesetting}[filename, in directory art]
	\includegraphics[width=\textwidth]{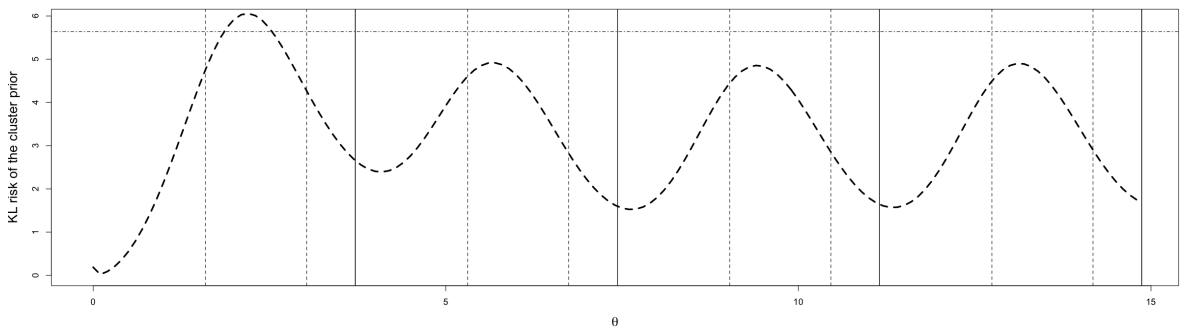}
	% note that files may not be rotated
	\caption{Plot of the univariate predictive KL risk $\rho(\theta,\hat p_{\sf C}[\eta,r]$ as $\theta$ varies over the $x$-axis. Here, $\eta=0.001$ and $r=0.225$. The horizontal line corresponds to the asymptotic minimax limit $\lambda^2_f(\eta)/(2r)$. The dotted vertical lines denotes the location of the non-origin support points of  $\pi_{\sf C}[\eta,r]$ with the bold lines marking each cluster boundary.}
	\label{fig2}
\end{figure}

Figure~\ref{fig2} shows the numerical evaluation of the predictive risk $\rho(\theta,\hat p_{\sf C}[\eta,r])$ of the cluster prior based Bayes \textit{prde}  when $\eta=0.001$ and $r=0.225$. Each cluster has size three. %The figure demonstrates the non-asymptotic behavior of the risk bounds used in the proof of Theorem 1.  
The maximum risk $\hat p_{\sf C}[\eta,r]$ crosses the  asymptotic theory limit but does not exceed by much. It  shows that the asymptotic analysis is fairly reflective in this non-asymptotic regime.  The risk function has its peak between $\mu_{11}$ and $\mu_{12}$ and is approximately periodic barring a few clusters near the origin. As the figure shows, the risk function is much smaller than the asymptotic limit of $\lambda_{f}^2/(2r)$ for all the points in $C_1$ barring its first point. As all points in $C_1$ are equally likely, this implies that the cluster prior is not least favorable. The following result make this observation rigorous by explicitly evaluating the first order asymptotic Bayes risk of the cluster prior. It establishes that when there are two or more points in each cluster (i.e. $r< r_0$)  the cluster prior is no longer least favorable. Its Bayes risk, however, has the same order of the minimax risk and will be at least 34\% of the minimax risk for any value of $r$. 

\bt 
	If $\eta_n=s_n/n \to 0$ as $n \to \infty$, then the multivariate cluster prior $\pi^{\sf C}_n[\eta_n,r]$ is not asymptotically least favorable for all $r < r_0$. As such, its Bayes risk satisfies:
	\[
	\lim_{n \to \infty} \; \bigg\{B(\pi_n^{\sf C}[\eta_n,r]) \big)\bigg\}\bigg/ R^*(\Theta_0[s_n])=\frac{1}{K_r} \bigg \{ 1+r \sum_{i=1}^{\infty} \big(1+r^{-1}-(1+4r)^{2i}\big)_+\bigg\}~,
	\]
	where, $K_r$ is defined in \eqref{eq.K}. Additionally, if $\eta_n \to 0$ and $s_n \to \infty$ as $n \to \infty$ then $\pi^{\sf C}_n[\eta_n,r]$ is asymptotically least favorable for all $r\geq r_0$.
	%$$\liminf_{n \to \infty}  B(\pi^{\sf C}_n)/ R^*(\Theta_0[s_n])=1.$$
\et
\hspace{1cm}\\[-9ex]
\section{Proof Layout} 
We provide a brief overview of the proof of our main result. Detailed proofs are provided in the supplement. The proof of Theorem 1 involves asymptotically upper bounding the risk  $\sup_{\theta \in \Theta_0[s_n]} \rho(\theta,\hat p_{\sf C})$ by $R^*(\Theta_0[s_n])$. Then, the asymptotic equality follows as the first term can not be smaller than the minimax risk by definition. Also,  note that due to the product structure of the prior, the multivariate maximal  risk can be evaluated based on the  risk of the univariate Bayes \textit{prde} $\hat p_{\sf C}[\eta_n,r]$ by using the following relation: 
\begin{align}\label{eq:mult}
\sup_{\theta \in \Theta_0[s_n]} \rho(\theta,\hat p_{\sf C})=n(1-\eta_n) \rho(0,\hat p_{\sf C}[\eta_n,r]) + n \eta_n \sup_{\theta \in \mathbb{R}\setminus 0}\rho(\theta,\hat p_{\sf C}[\eta_n,r])~.
\end{align}
Asymptotic evaluation of the two expressions on the right above is done by using the risk decomposition lemma 2.1 of MJ17. It reduces the calculation for the univariate predictive risk to finding expectation of functionals involving standard normal random variable $Z$ as 
\begin{align}\label{eq:risk-decom}
&\rho(\theta,\hat{p}_{\sf C}[\eta_n,r])=\frac{\theta^2}{2r}-\Exp\log N_{\theta,v}(Z)+\Exp\log D_\theta(Z), \text{ where,} \\
& N_{\theta,v}(Z) = 1 + \sum_{i \in \mathbb{Z}\setminus 0} \frac{q_{i}}{K}\sum_{j=1}^K N_{ij}(\theta, Z;v)
\text{ and }
D_\theta(Z)=N_{\theta,1}(Z)~.\nonumber
\end{align}
Here, $q_{i}=(1-\eta_n)^{-1}P(C_i)$ with $P(C_i)$ being the mass of cluster $C_i$ in $\pi_{\sf C}[\eta_n,r]$;  thus $q_i =2^{-1}\exp(-|i|\lambda_{e,n}^2/2)$ with $\lambda_{e,n} =(2 \log \eta_n^{-1})^{-1}$ and $\lambda_{f,n}=v^{1/2} \lambda_{e,n} $; $N_{ij}$ is the contribution to the risk of the $j$th support point $\mu_{ij}(\eta_n,r)$ within the $i$th cluster.
 
%For ease of presentation, we do not make the dependence of $q_i, \mu_{ij}$ on $\eta_n$ explicit in the notations. 
The risk contributions $N_{ij}$ are exponents of quadratic forms in $\mu_{ij}$, viz,  
%\[
$N_{ij}(\theta, Z;v)=\exp \{v^{-1/2}{\mu_{ij}Z}+ v^{-1}\mu_{ij}\theta -(2v)^{-1}{\mu_{ij}^2}\}$. The risk at the origin is well-controlled for this cluster prior based \textit{prde} (lemma 1  of supplement) and so, based on
%\]
\eqref{eq:mult}, it is suffices to bound  $\sup_{\theta} \rho(\theta,\hat{p}_{\sf C}[\eta_n,r])$ by $ \lambda_{f,n}^2/(2r)$ to arrive at the desired result. This involves tracing two fundamentally different risk  phenomena depending on the location of $\theta$ (a) $\theta \in C_{\pm 1}$ (b) $\theta \notin C_{\pm 1}$. 
In the former case,  $\Exp\log D_\theta(Z)=O(\lambda_{f,n})$ (by lemma 3 of the supplement) and thus the contribution of the third term on the right of \eqref{eq:risk-decom} is not significant. 
Also, $\Exp\log N_{\theta,v}(Z)=O(\lambda_{f,n})$ for $|\theta|\leq \lambda_{f,n}$ and so, asymptotically $\rho( \theta,\hat{p}_{\sf C}[\eta_n,r])$ initially increases quadratically in $\theta$ and $\rho( \lambda_{f,n},\hat{p}_{\sf C}[\eta_n,r])=\lambda_{f,n}^2/(2r)(1+o(1))$. However, if $|\theta| \in C_1 \setminus [0, \lambda_{f,n}] $, then  $\Exp\log N_{\theta,v}(Z)$  is significantly large and controls the predictive risk below the desired asymptotic limit (see lemma~4 of supplement). 

If $\theta \in C_i$ for any $|i|>1$, then the risk phenomenon is quite different than the origin adjoining clusters. Now,  $\Exp\log D_\theta(Z)$ is significantly positive. However, an important ingredient of the proof is that its magnitude can be asymptotically well controlled by considering only atoms in $C_i$ or the nearest atom in $C_{i-1}$. Lemma~3 in the supplementary material establishes that for $\theta \in C_i$ with $|i|>1$, 
$\Exp\log D_\theta(Z)\leq \{\Exp\log D_{i.}(Z)\}_+ + o(\lambda^2_{f,n} ) \text{ where } D_{i.}(Z)=N_{i-1,K}(\theta,Z;1)+\sum_{j=1}^K N_{ij}(\theta,Z;1).$ 
Next, use the naive bound $\Exp\log N_{\theta,v}(Z)\geq \Exp\log N_{i.}(Z))$ where $N_{i.}=N_{i-1,K}(\theta,Z;v)+\sum_{j=1}^K N_{ij}(\theta,Z;v)$. Now, plugging these two bounds in \eqref{eq:risk-decom} we get the desired upper bound  (see lemma~4 of the supplement).

\section{Simulations}
\begin{table}[b]
	\centering
	\caption{Numerical evaluation of the maximum risk for the different univariate predictive
		density estimates as the degree of sparsity ($\eta$) and predictive difficulty $r$ varies.
		The asymptotic minimax risk is reported in `A-Theory' and the subsequent columns report
		the maximum risk of the estimators as quotients of `A-Theory' values. }
	\scalebox{0.74}{
	\begin{tabular}{|c|c|cccccccc|}
		\toprule
		\textbf{Sparsity} & \textbf{r} & \;\textbf{A-Theory}\; & \;\textbf{Plugin}\; & \;\textbf{Thresh}\; & \;\textbf{E-Grid}\; & \;\textbf{P-Grid}\; & \;\textbf{Bi-Grid}\; & \;\textbf{SUS}\; & \;\textbf{Clustered}\; \\
		\midrule
		& 1     & 2.3026 & 1.0841 & 0.7057 & 0.6236 & 0.7366 & 0.7366 & 0.9090 & 0.7629 \\
		\cmidrule{2-10}          & 0.5   & 3.0701 & 1.6023 & 0.8822 & 0.8031 & 0.8832 & 0.8832 & 1.0135 & 1.2036 \\
		\cmidrule{2-10}    0.01  & 0.25  & 3.6841 & 2.6310 & 0.9235 & 1.2718 & 1.0398 & 1.0079 & 1.1383 & 1.0932 \\
		\cmidrule{2-10}          & 0.1   & 4.1865 & 5.6949 & 1.1074 & 2.6198 & 1.2304 & 1.2239 & 1.2677 & 1.3507 \\
		\midrule
		& 1     & 5.7565 & 1.1371 & 0.7332 & 0.7407 & 0.7277 & 0.7277 & 0.8665 & 0.7287 \\
		\cmidrule{2-10}          & 0.5   & 7.6753 & 1.6960 & 0.8522 & 0.9543 & 0.8486 & 0.8486 & 0.9599 & 1.0874 \\
		\cmidrule{2-10}    0.00001 & 0.25  & 9.2103 & 2.8120 & 0.9125 & 1.4146 & 0.9781 & 0.9464 & 1.0328 & 1.0376 \\
		\cmidrule{2-10}          & 0.1   & 10.4663 & 6.1542 & 1.0395 & 2.7946 & 1.1049 & 1.0710 & 1.1182 & 1.0932 \\
		\midrule
		& 1     & 11.5129 & 1.2390 & 0.7958 & 0.8357 & 0.7891 & 0.7891 & 0.8765 & 0.7910 \\
		\cmidrule{2-10}    1.00E-10 & 0.5   & 15.3506 & 1.8540 & 0.8810 & 1.0488 & 0.8734 & 0.8734 & 0.9337 & 1.1080 \\
		\cmidrule{2-10}          & 0.25  & 18.4207 & 3.0835 & 0.9451 & 1.5092 & 0.9855 & 0.9629 & 0.9945 & 1.0128 \\
		\cmidrule{2-10}          & 0.1   & 20.9326 & 6.7701 & 1.0191 & 2.8958 & 1.1008 & 1.0138 & 1.0611 & 1.0233 \\
		\bottomrule
	\end{tabular}%
	}\label{tab-1}	
\end{table}%

We introspect the performance of the aforementioned \textit{prde}s across different sparsity regimes. The product structure of our estimation framework allows us to concentrate on the maximal risk of the corresponding univariate \textit{prde}s. In table~\ref{tab-1}, we report the maximum risk of our proposed clustered prior based Bayes (CB) \textit{prde} (in last column) as the degree of sparsity $\eta$ and predictive difficulty $r$ varies. The performance of the six following competing methods 
(a) hard thresholding based plugin estimator (b) thresholding based risk diversified \textit{prdre} of MJ15 and
Bayes \textit{prdes} based on (c)   $\pi_{\sf EG}$ prior of J94 (d) $\pi_{\sf PG}$ prior of MJ17 (e) $\pi_{\sf BG}$ prior of MJ17 (f) spike and uniform slab (SUS) prior, are respectively reported in columns 4 to 9 in table~\ref{tab-1}. Across all regimes the maximum risk of CB-\textit{prde} is  reasonably close to the order of the minimax risk prescribed by the asymptotic theory; for large $r$ values the maximum risk is actually lower than the asymptotic theory prescribed minimax value whereas it is little higher for lower $r$ values, particularly at moderate sparsity.   
For lower $r$ values, CB-\textit{prde} is substantially better than that the plugin or grid prior based \text{prde}s. Overall, CB-\textit{prde} has similar performance to that of the risk diversified \textit{prde}s of MJ15 and MJ17, both of which are asymptotically minimax optimal for all $r$. 

%There is no need to have clusters with $K$ points for $C_i>1$. Having only a point with $\eta^i_n$ will be sufficient. 

\section*{Supplementary Materials and Acknowledgement}
Detailed proofs of the results stated in Section 1 are provided in the supplementary materials. 
GM is indebted to Professor Iain Johnstone for numerous stimulating discussions which led to many of the ideas in this paper. The research here was partially supported by NSF DMS-1811866.\\[-5ex]

%\clearpage 
\bibliography{pred-inf,append}

\end{document}

\title[Sparse Minimax Optimality of Bayes PDE from Clustered Discrete Priors]{\texttt{Supplement to} ``Sparse Minimax Optimality of Bayes Predictive Density Estimates from Clustered Discrete Priors"}
\date{}
%\author{Ujan Gangopadhyay}
%\address{Ujan Gangopadhyay, \ Department of Mathematics, \ University of Southern California, \ Los Angeles, CA, USA.}
%\email{ujan.gangopadhyay@usc.edu}
%\author{Gourab Mukherjee}
%\address{Gourab Mukherjee, \ Marshall School of Business, \ University of Southern California, \ Los Angeles, CA, USA.}
%\email{gourab@usc.edu}

\begin{abstract}
This supplement contains detailed proofs of the results in the paper. We first provide some background and preliminary results on predictive KL risks of Bayes predictive density estimators. Thereafter, we provide detailed proofs of the theorems in the main paper with separate analysis for the sub-critical case $r<r_0$ and the super-critical case $r\geq r_0$. We also explain why the choice of $r_0=1/2$ is optimal.
\end{abstract}

\maketitle

\bibliographystyle{plainnat}

%\begin{keywords}
%Address; Appendix; Figure; Length; Reference; Style; Summary; Table.
%\end{keywords}

\section{Background and Preliminaries}

For the technical proofs without loss of generality assume $v_x=1$. So, $r=v_x/v_y=v_x$. Recall $v=(1+r^{-1})^{-1}$ and $\eta_n=s_n/n$. As demonstrated in equation  (3) of the main paper, the multivariate maximal risk of the Bayes predictive density estimate (\textit{prde}) from the cluster prior can be evaluated by studying the  predictive risk of the univariate Bayes \textit{prde} $\hat p_{\sf C}[\eta_n,r]$ based on the univariate cluster prior $\pi_{\sf C}[\eta_n,r]$:  
\beq\sup_{\theta \in \Theta_0[s_n]} \rho(\theta,\hat p_{\sf C})=n(1-\eta_n) \rho(0,\hat p_{\sf C}[\eta_n,r]) + n \eta_n \sup_{\theta \in \mathbb{R}\setminus 0}\rho(\theta,\hat p_{\sf C}[\eta_n,r])~.\label{Eq: mvt}\eeq 
Henceforth, unless we explicitly mention, we would concentrate on univariate Bayes predictors and their risk functions. 
Recall, in the multivariate set-up we consider asymptotically sparse regimes, where $\eta_n\to 0$ as $n\to\infty$. 
Hereon, for notational convenience we write $\eta$ instead of $\eta_n$ keeping the dependence on $n$ implicit. Recall,
\[
\lame:=\sqrt{2\log\eta^{-1}},\quad\mbox{and}\quad
\lamf:=\sqrt{2\log\eta^{-v}}.
\]
Recall from equation (2) of the main paper, the univariate clustered discrete prior $\pi_{\sf C}[\eta,r]$ is the following:  
\[
\pi_{\sf C}[\eta,r]=(1-\eta)\delta_0+\frac{1-\eta}{2} \sum_{j=1}^{\infty} \eta^{j} \big\{C_j (\eta,r) + C_{-j}(\eta,r)\big\}.
\]
The point-masses in cluster $C_j$ are denoted by $\{\mu_{jk}: k=1,\ldots,K\}$ where the common cluster size $K$ is
\[
K:=K(r)=1+   \big \lceil \log(1+r^{-1})/(2\log(1+4r))\big \rceil  \cdot 1\{r< r_0\},
\]
where, $r_0=1/2$. Further, recall that 
\[
C_j(\eta,r)= \frac{1}{K} \sum_{k=1}^{K}  \delta_{\mu_{jk}} \text{ for } j \in \mathbb{Z}\setminus \{0\},
\]
where $\mu_{1k}= \lamf (1+4 r)^{k-1} \wedge \lambda_e$ for $1\leq k \leq K$, $\mu_{jk}=(j-1)\mu_{1K}+\mu_{1k}$ for $j\geq 2$, and $\mu_{jk}=\mu_{-jk}$ for $j<0$. So for $r\geq r_0$, that is, when $K=1$, the clustered discrete prior only has point-masses $\LP j\lamf: j\in\ZZ\RP$.

By Lemma 2.1 of \citet{mukherjee2017minimax} the predictive KL risk of the univariate cluster prior is given by: 
\beq
\rho(\theta,\hat{p}_{\sf C}[\eta,r])=\frac{\theta^2}{2r}-\Exp\log N_{\theta,v}(Z)+\Exp\log D_\theta(Z)
\label{Eq: decomp}
\eeq
where $Z$ is a standard normal random variable, and
\begin{align*}
N_{\theta,v}(Z) & = 1 + \frac{1}{2K}\sum_{j\in\ZZ\bs\LP 0\RP}\sum_{k=1}^K\exp\LP\frac{\mu_{jk}Z}{\sqrt{v}}+\frac{\mu_{jk}\theta}{v}-\frac{\mu_{jk}^2}{2v}-|j|\frac{\lame^2}{2}\RP,\mbox{ and}\\
D_{\theta}(Z) & = 1 + \frac{1}{2K}\sum_{j\in\ZZ\bs\LP 0\RP}\sum_{k=1}^K\exp\LP\mu_{jk}(Z+\theta)-\frac{\mu_{jk}^2}{2}-|j|\frac{\lame^2}{2}\RP.
\end{align*}

\section{Proof of Theorem 1}

We first present the proof for $r<r_0$ because the proof is more intricate compared to the case when $r\geq r_0$. In the latter case, by definition $K=1$, and the proof is is comparatively easier. It uses parts of the proof techniques used for $r<r_0$ case but also involves some fundamentally different attributes. Hence, it is presented afterwards where we also explain the choice of $r_0=1/2$.

\subsection{Notations}\label{subsec: notation}

For convenience of notation, we shall write the support points of the clustered discrete prior as $\LP\mu_p:p\in\ZZ\RP$. The identification is made as follows. Let $\mu_0=0$, and for $p>0$ identify $\mu_p$ in the new notation with $\mu_{j_p k_p}$ in the original notation where $j_p,k_p$ are the unique positive integers such that $p=(j_p-1)K+k_p$ with $k_p\leq K$. For $p<0$ let $\mu_p=-\mu_{-p}$. So essentially $\mu_p$ is the $k_p$th point in the $p$th cluster. Let $j_0=0$ and $j_{-p}=j_p$ for $p<0$. Let $c_0=1$ and $c_p=(2K)^{-1}$ for $p\neq 0$. With these new notations can write 
\[
D_\theta(Z) = \sum_{p\in\ZZ} D_{\theta p}(Z),\mbox{ and }
N_\theta(Z) = \sum_{p\in\ZZ} N_{\theta p}(Z)
\]
where
\begin{align*}
D_{\theta p}(Z) &:= c_p\exp\LP\mu_p Z+\mu_p\theta-\frac{1}{2}\mu_p^2-j_p\frac{\lame^2}{2}\RP,\mbox{ and}\\
N_{\theta p}(Z) &:= c_p\exp\LP\frac{\mu_p Z}{\sqrt{v}}+\frac{\mu_p\theta}{v}-\frac{\mu_p^2}{2v}-j_p\frac{\lame^2}{2}\RP.
\end{align*}
The above notations will be used for all $r\in(0,\infty)$. But now we define two indexes $\ld(\theta)$ and $\ln(\theta)$ for all $\theta>0$ specifically for $r<r_0$. If $\theta\in[j\lame,(j+1)\lame)$, then let $\ld(\theta):=jK$. So $\ld(\theta)$ is the number of support points in the cluster prior between $[0,j\lame]$. This is the index of the atom $\mu_p$ such that $\Exp D_{0p}(Z)$ is maximized. Note that $j\lame=\mu_{jK}=\mu_{\ld(\theta)}$. Now we define the index $\ln(\theta)$ which is the index of the atom $\mu_p$ such that $\Exp N_{0 p}(Z)$ is maximized. More precisely, $\ln(\theta)$ is defined as follows:   

\bei 
\ii If $\theta\in[j\lame,j\lame+\lamf]$, then let $\ln(\theta):=jK$. Note that, in this case $\mu_\ln=\mu_{jK}=j\lame$.

\ii If $\theta\in(j\lame+\lamf(1+4r)^k,\min\{j\lame+\lamf(1+4r)^k(1+2r),(j+1)\lame\}]$ for $0\leq k<K$, then let $\ln(\theta):=jK+k+1$. Note that, in this case $\mu_\ln=\mu_{jK+k+1}=j\lame+\lamf(1+4r)^k$.

\ii If $\theta\in(j\lame+\lamf(1+4r)^k(1+2r),\min\{j\lame+\lamf(1+4r)^{k+1},(j+1)\lame\}]$ for some $0\leq k<K$, then let $\ln(\theta):=jK+k+2$. Note that, $\mu_\ln=\mu_{jK+k+2}=\min\{j\lame+\lamf(1+4r)^{k+1},(j+1)\lame\}$.
\ee

\subsection{Risk at origin}

The risk at the origin for our cluster prior based Bayes \textit{prde} is asymptotically much smaller than the risk for the thresholding based risk diversified \textit{prde} of \citet{Mukherjee-15}. As such, comparing equation (51) in the aforementioned paper with the following result, it follows that any thresholding based minimax optimal \textit{prde} will have much higher risk at the origin than the cluster prior based Bayes \textit{prde}. The Bayes \textit{prde}s based on grid and bi-grids priors such as the $\pi_{\sf EG}$ prior of \citet{Johnstone94a} and $\pi_{\sf PG}$ and $\pi_{\sf PG}$ priors of \citet{mukherjee2017minimax} have similar risk to the cluster prior based Bayes \textit{prde} at the origin. 

\begin{lemma} For any fixed $r \in (0,\infty)$, $\rho(0,\hat p_{\sf C}[\eta,r])\leq \eta(1+o(1))$ as $\eta\to 0$.\label{Lemma: Origin}\end{lemma}

\bpf By definition $N_{\theta, v}(Z)\geq 1$ for all $Z$. Using \eqref{Eq: decomp}, we have
\[
\rho\Lp 0,\hat p_{\sf C}\Rp = -\Exp\log N_{\theta,v}(Z) + \Exp\log D_\theta(Z) \leq \Exp\log D_\theta(Z).
\]
Note that, for $p\neq 0$, $\Exp D_{0 p}(Z)=(2K)^{-1}\eta^{j_p}$. Summing over all $p\neq 0$ and using $D_0=1$ along with the inequality $\log(1+x)\leq x$ for $x\geq 0$ we get
\[
\Exp\log D_{0}(Z)\leq \sum_{p\neq 0}\Exp D_{0p}(Z)=\sum_{p\neq 0}(2K)^{-1}\eta^{j_p}=\sum_{j=1}^\infty\eta^j=\frac{\eta}{1-\eta}.
\]
This completes the proof.
\epf

\subsection{Risk bounds at the non-origin parametric points}

Next, we concentrate on the risk at the non-origin points. Our goal is to establish
\beq
\sup_{\theta\in\RR\bs\LP 0\RP}\rho\Lp\theta,\hat p_{\sf C}[\eta,r]\Rp\leq \frac{\lamf^2}{2r}(1+o(1))\quad\mbox{as $\lamf\to\infty$}.
\label{Eq: nonorigin}\eeq
This along with \eqref{Eq: mvt} and the above result about the risk bound at the origin will imply that the multivariate maximum risk obeys
\begin{align*}
\sup_{\theta\in\Theta_0[s_n]}\rho(\theta,\hat p_{\sf C}[\eta_n,r]) & \leq
-n\eta_n(1-\eta_n)^{-1}+n\eta_n\frac{\lamf^2}{2r}(1+o(1))\\ & = 
n\eta_n\log\eta_n^{-1}(1+r)^{-1}(1+o(1))
\end{align*}
which would establish the result in Theorem 1. By symmetry, it would be enough to prove the bound in \eqref{Eq: nonorigin} for positive $\theta$. Hence, hereon in this subsection we only consider $\theta>0$. In this case the contribution of $D_{\theta i}$s for $i<0$ are expected to be negligible. This is formalized in the following result.

\begin{lemma} For any $r\in(0,\infty)$ and any fixed $\theta>0$ we have,
\[
\Exp\log D_\theta(Z) = \Exp\log\Lp 1+\sum_{i=1}^\infty D_{\theta i}(Z)\Rp + o(1)\quad\mbox{as $\lamf\to\infty$.}
\]
\label{Lemma: Negative}\end{lemma}

\bpf Using the the inequality $\log(1+x+y)\leq x+\log(1+y)$ for nonnegative $x,y$ we get,
\[
\Exp\log D_\theta(Z) \leq \sum_{i<0}\Exp D_{\theta i}(Z) 
+ \Exp\log\Lp 1+\sum_{i=1}^\infty D_{\theta i}(Z)\Rp.
\]
Using definition of $j_i$ and $D_{\theta i}$ and the fact that $\mu_i<0,\theta>0$ we get,
\[
\Exp D_{\theta i}(Z) = \Exp\frac{1}{2}\exp\LP\mu_i(Z+\theta)-\frac{\mu_i^2}{2}-j_{i}\frac{\lame^2}{2}\RP=\frac{1}{2}\exp\LP\mu_i\theta-j_{i}\frac{\lame^2}{2}\RP\leq \frac{1}{2}e^{-j_{i}\frac{\lame^2}{2}}.
\]
As $i$ runs from $0$ to $-\infty$, $j_{i}$ goes from $1$ to $\infty$ with each term repeating $K$ times. Hence, summing over $i<0$ we get, $\sum_{i<0}\Exp D_{\theta i}=o(1)$ as $\lamf\to\infty$. This completes the proof.
\epf 

We first provide an upper bound on $\Exp\log D_\theta(Z)$, which would be substituted in equation \eqref{Eq: decomp} to get the required upper bound. The following result is crucial as it shows that the infinite sum in the expression of $D_\theta(Z)$ can be asymptotically reduced as a contribution from a single dominant term. This reduction greatly helps in tracking the risk of the cluster prior and is pivotal in the proof of Theorem 1.

\begin{lemma} For $r<r_0$ and any fixed $\theta>0$ we have, 
\[
\Exp\log D_\theta(Z) = \Exp\log D_{\theta\ld}(Z)+ O(\lamf)\quad\mbox{as $\lamf\to\infty$}.
\]
\label{Lemma: D}\end{lemma} 

\bpf By virtue of the previous lemma, we can consider only contributions from $D_{\theta i}(Z)$s with $i>0$. We suppress the dependence of $D_{\theta i}(Z)$s on $\theta$ and $Z$ and simply write $D_i$. Similarly $D_\theta(Z)$ is written only as $D_\theta$. Note that, $\ld\geq 0$ because $\theta>0$. First we get an upper bound on $\Exp\log D_\theta$ by separating the contribution from $\mu_\ld$ as follows
\beq
\Exp\log\Lp 1+\sum_{i=1}^\infty D_{\theta i}(Z)\Rp\leq
\Exp\log D_\ld + 
\Exp\log\Lp 1+\sum_{i=\ld+1}^\infty\frac{D_i}{D_\ld}\Rp + 
\Exp\log\Lp 1+\sum_{i=0}^{\ld-1}\frac{D_i}{D_\ld}\Rp.
\label{Eq: aux2}
\eeq
In the right hand side of the above equation, the second term compares contribution of $\mu_\ld$ with that of the succeeding terms. We split it further by separating out the contribution of points in the next cluster from the rest in the following manner
\beq
     \Exp\log\Lp 1+\sum_{i=\ld+1}^\infty\frac{D_i}{D_\ld}\Rp 
\leq \sum_{i=\ld + 1}^{\ld+K}\Exp\log\Lp 1+\frac{D_i}{D_\ld}\Rp+
\sum_{i=\ld+K}^\infty\Exp\frac{D_{i+1}}{D_i}.
\label{Eq: aux1}
\eeq 
Note that by definition of $\ld$, $\mu_\ld<\theta\leq\mu_{\ld+K}$. Take $i$ such that $\ld+1\leq i\leq \ld+K$. Let $d_i:=\mu_i-\mu_{\ld}$. Using the inequality $\log(1+x)\leq\log 2+(\log x)_+$ we get  
\begin{multline}
\Exp\log\Lp 1+\frac{D_i}{D_\ld}\Rp
= \Exp\log\Lp 1+\exp\LP d_{i}\Lp Z+\theta-\mu_\ld-\frac{d_{i}}{2}\Rp-\frac{\lame^2}{2}\RP\Rp\\
\leq \Exp\log\Lp 1+\exp\LP d_{i} Z -\frac{1}{2}(\lame-d_{i})^2\RP\Rp 
\leq \log 2 + \Exp(d_{il} Z-(\lame-d_{i})^2/2)_+=O(\lamf).
\label{Eq: aux3}
\end{multline}
Summing over $\ld+1\leq i\leq \ld+K$ we get the fist term in the right-hand side of equation \eqref{Eq: aux1} is $O(\lamf)$. Now, to consider the second term in the right-hand side of \eqref{Eq: aux1}, take $i\geq\ld+K$. Using $\Exp((\mu_{i+1}-\mu_{i})Z)=(\mu_{i+1}-\mu_{i})^2/2$ we get
\beq
\Exp\frac{D_{i+1}}{D_i}\leq\Exp\exp\LP(\mu_{i+1}-\mu_i)\Lp Z+\theta-\frac{\mu_{i+1}-\mu_i}{2}\Rp\RP
=\exp\LP(\mu_{i+1}-\mu_i)(\theta-\mu_i)\RP.
\label{Eq: aux7}\eeq
Since $i$ runs from $\ld+K$ to $\infty$ and $\theta\leq\mu_{\ld+K}$, we get the second term in the right-hand side of equation \eqref{Eq: aux1} is $O(1)$. Hence, the second term in the right-hand side of equation \eqref{Eq: aux2} is $O(\lamf)$. Now we consider the third term in the right-hand side of equation \eqref{Eq: aux2} is also $O(\lamf)$. We split the sum as
\begin{multline}
\Exp\log\Lp 1+\sum_{i=0}^{\ld-1}\frac{D_i}{D_\ld}\Rp \leq 
\sum_{i=\ld-K+1}^{\ld-1}\Exp\log\Lp 1+\frac{D_i}{D_\ld}\Rp\\ +
\Exp\log\Lp 1+\frac{D_{\ld-K}}{D_\ld}\Rp +
\sum_{i=0}^{\ld-K-1}\Exp\frac{D_i}{D_{\ld-K}}.
\label{Eq: aux4}
\end{multline}
To consider the first term in the right-hand side above, take $\ld-K+1\leq i\leq \ld-1$. Then $\theta\geq\mu_\ld\geq(\mu_\ld+\mu_i)/2$ and because of the structure of the atoms in the clusters, $\theta-(\mu_\ld+\mu_i)/2=\Theta(\lamf)$. Note that, $i$ and $\ld$ belong to the same cluster. Using symmetry of the distribution of $Z$ we get
\begin{multline*}
\Exp\log\Lp 1+\frac{D_i}{D_\ld}\Rp=
\Exp\log\Lp 1+\exp\Lp(\mu_i-\mu_\ld)\Lp Z+\theta-\frac{\mu_\ld+\mu_i}{2}\Rp\Rp\Rp\\=
\Exp\log\Lp 1+\exp\Lp(\mu_\ld-\mu_i)\Lp Z-\theta+\frac{\mu_\ld+\mu_i}{2}\Rp\Rp\Rp = O(1). 
\end{multline*}
Summing over $\ld-K+1\leq i\leq \ld-1$ we get the first term in the right-hand side of equation \eqref{Eq: aux4} is $O(1)$. Now consider the second term in the right-hand side of equation \eqref{Eq: aux4}. As before
\[
\Exp\log\Lp 1+\frac{D_{\ld-K}}{D_\ld}\Rp 
=\Exp\log\Lp 1+\exp\Lp(\mu_{\ld}-\mu_{\ld-K})\Lp Z+\frac{\mu_\ld+\mu_{\ld-K}}{2}-\theta\Rp+\frac{\lame^2}{2}\Rp\Rp.
\]
Note that, $\mu_\ld-\mu_{\ld-K}=\lame$ and $(\mu_\ld+\mu_\ld-K)/2-\theta\leq-\lame/2$. Therefore,
\begin{align*}
\Exp\log\Lp 1+\frac{D_{\ld-K}}{D_\ld}\Rp & = 
\Exp\log\Lp 1+\exp\Lp\lame\Lp Z-\theta+\frac{\mu_\ld+\mu_{\ld-K}}{2}\Rp+\frac{\lame^2}{2}\Rp\Rp\\ & \leq
\Exp\log\Lp 1+\exp\Lp\lame Z\Rp\Rp\leq
\log 2+\lame\Exp Z_+=O(\lamf).
\end{align*}
Finally, for each $0\leq i<\ld-K$ define $b_i=\lfloor(\ld-K-i)/K\rfloor$. Then  
\[
\Exp\log\Lp 1+\frac{D_i}{D_{\ld-K}}\Rp\leq
\Exp\frac{D_i}{D_{\ld-K}}\leq
2K\exp\Lp(\mu_{\ld-K}-\mu_i)(\mu_{\ld-K}-\theta)+b_i \frac{\lame^2}{2}\Rp.
\]
Note that, $\theta-\mu_{\ld-K}\geq\lame$ and $\mu_{\ld-K}-\mu_i\geq b_i\lame$. Thus, we have
\[
\Lp\mu_{\ld-K}-\mu_i\Rp\Lp\mu_{\ld-K}-\theta\Rp+b_i \frac{\lame^2}{2}\leq-b_i\lame^2+b_i\frac{\lame^2}{2}=-b_i\frac{\lame^2}{2}.
\]
Summing over all $0\leq i\leq\ld-K$ we arrive at the following bound
\[
\sum_{0\leq i\leq\ld-K}\Exp\frac{D_i}{D_{\ld-K}}\leq K\sum_{p=0}^{b_0}e^{-p\lame^2/2}=O(1).
\]
This shows that the third term in the right-hand side of \eqref{Eq: aux4} is $O(1)$. Thus we have proved that the second and third term in the right-hand side of \eqref{Eq: aux2} are $O(\lamf)$ as $\lamf\to\infty$. This completes the proof.
\epf

The previous lemma essentially shows that to get an upper bound on $\Exp\log D_\theta(Z)$ it is enough to consider only $D_{\theta\ld}(Z)$ because asymptotically the contribution of the other terms are negligible. To prove \eqref{Eq: nonorigin} using \eqref{Eq: decomp} we need a lower bound on $\Exp\log N_{\theta,v}(Z)$, which we get by the straightforward inequality $\Exp\log N_\theta(Z) \geq \Exp\log N_{\theta\ln}(Z)$. Of course the novelty is in choice of $\ln(\theta)$ and in the next result we see that these bounds are enough to prove \eqref{Eq: nonorigin}. 

\begin{lemma} For $r<r_0$ and for any $\theta>0$, with $\ln(\theta)$, $\ld(\theta)$, $N_{\theta\ln}$, $D_{\theta\ld}$ defined in Subsection~\ref{subsec: notation} we have 
\[
\frac{\theta^2}{2r} - \Exp\log N_{\theta\ln}(Z) + \Exp\log D_{\theta\ld}(Z) \leq \frac{\lamf^2}{2r}\Lp 1+o(1)\Rp\mbox{ as }\lamf\to\infty.
\]
\label{Lemma: DNCancel}\end{lemma}

\bpf For convenience we write $\ld(\theta)$ and $\ln(\theta)$ as $\ld$ and $\ln$ respectively. Note that from the definition of $\ld$ it follows $\mu_\ld\leq\theta\leq\mu_{\ld+K}$. Let 
\[
A_\theta:=\mu_\ld\theta-\frac{\mu_\ld^2}{2}-j_\ld\frac{\lame^2}{2},
\quad\mbox{and}\quad 
B_\theta:=\frac{\mu_\ln\theta}{v}-\frac{\mu_\ln^2}{2v}-j_\ln\frac{\lame^2}{2}.
\]
From definitions it follows $D_\ld = c_\ld\exp(\mu_\ld Z + A_\theta)$ and $N_\ln = c_\ln\exp(\mu_\ln Z+B_\theta)$. Hence,
\beq
-\Exp\log N_{\theta \ln}(Z) \leq - B_\theta + O(1)\text{ as $\lamf\to\infty$.}
\label{Eq: aux5}
\eeq
Using $\theta\geq\mu_\ld=j_\ld\lame$ we get
\[
A_\theta=j_\ld\theta-\frac{j_\ld^2\lame^2}{2}-j_\ld\frac{\lame^2}{2}\geq\Lp j_\ld^2-j_\ld\Rp\frac{\lame^2}{2}\geq 0.
\]
Using this, we derive the upper bound
\begin{multline*}
\Exp\log(1+D_{\theta\ld}(Z))=
\Exp\log(1+c_\lambda\exp(\mu_\ld Z + A_\theta))\\=
A_\theta+\Exp\log(c_\lambda+\exp(\mu_\ld Z-A_\theta))\leq
A_\theta+\Exp(\mu_\ld Z - A_\theta)_+ + O(1).
\end{multline*}
Since $\mu_\ld=j_\ld\lame$ and $A_\theta\geq(j_\ld^2-j_\ld)\lame^2/2$, we see that $\Exp(\mu_\ld Z-A_\theta)_+=O(\lamf)$. Hence,
\beq
\Exp\log(1+D_{\theta\ld}(Z))\leq A_\theta+O(\lamf).
\label{Eq: aux6}
\eeq
Combining equations \eqref{Eq: aux5} and \eqref{Eq: aux6} we get 
\[
\frac{\theta^2}{2r} - \Exp\log N_{\theta\ln}(Z) + \Exp\log D_{\theta\ld}(Z) \leq \frac{\theta^2}{2r}-A_\theta+B_\theta+O(\lamf).
\]
We will show that $\theta^2/(2r)+A_\theta-B_\theta\leq\lamf^2/(2r)$. First consider the case $\ln=\ld$ which means $\theta\in[j_\ld\lame,j_\ld\lame+\lamf]$. Observe in this case
\[
\frac{\lamf^2}{2r}-\frac{\theta^2}{2r}-A_\theta+B_\theta=\frac{\lamf^2}{2r}-\frac{\theta^2}{2r}+\frac{\mu_\ld\theta}{r}-\frac{\mu_\ld^2}{2r}=\frac{1}{2r}\Lp\lamf^2-(\mu_\ld-\theta)^2\Rp\geq 0.
\]
Now consider the case $\ln\neq\ld$, that is, $\ld+1\leq \ln\leq\ld+K$. In this case
\beq
\frac{\lamf^2}{2r}-\frac{\theta^2}{2r}-A_\theta+B_\theta=-\frac{\lamf^2}{2}-\frac{\theta^2}{2r}+\frac{\mu_\ln\theta}{v}-\frac{\mu_\ln^2}{2v}-\mu_\ld\theta+\frac{\mu_\ld^2}{2}.
\label{Eq: Quadratic}\eeq
This is a quadratic in $\theta$ if we fix values of $\ln$ and $\ld$. The roots are 
\begin{align*}
&\alpha_\ln:=\mu_\ln+r(\mu_\ln-\mu_\ld)-r\Lp\frac{(\mu_\ln-\mu_\ld)^2}{v}-\frac{\lamf^2}{r}\Rp^{1/2},\mbox{ and}\\
&\beta_\ln:=\mu_\ln+r(\mu_\ln-\mu_\ld)+r\Lp\frac{(\mu_\ln-\mu_\ld)^2}{v}-\frac{\lamf^2}{r}\Rp^{1/2}.
\end{align*}
Therefore, the quadratic in \eqref{Eq: Quadratic} is nonnegative in  the interval $[\alpha_\ln,\beta_\ln]$. So we need to check that for all the values of $\theta$ for which $\ln(\theta)=p$ lies in the interval $[\alpha_p,\beta_p]$. Because of the periodicity of the clusters, we can only consider the first cluster, so that $\ld=0$. In this case $\ln$ runs from $1$ to $K$. Using $\mu_0=0$ and $\mu_1=\lamf$ we get $\alpha_1=\lamf$ and $\beta_1=\lamf(1+2r)$. If $\beta_1>\mu_2$, then by definition $\ln$ is precisely $1$ and we are done. If $\beta_1<\mu_2\leq\lamf(1+4r)^2$, then $\mu_2=\lame$. After some calculations we see that $\alpha_2\leq\lamf(1+2r)$ with equality precisely when $\lame=\lamf(1+4r)^2$. In this case also, we see that the interval of $\theta$ for which $\ln(\theta)=2$ is contained inside the interval $[\alpha_2,\beta_2]$. In general for larger values of $K$, or equivalently, smaller values of $r$, it can be checked that for all $2\leq k\leq K$, $\alpha_k\leq\mu_k/(1+2r)$ and $\beta_k\geq\mu_{k}(1+2r)$. This shows the intervals $[\alpha_k,\beta_k]$ covers the set of $\theta$s for which $\ln(\theta)=k$ for all $k$. Therefore, the expression in \eqref{Eq: Quadratic} is nonnegative for all $\theta$, as required.  
\epf 

\subsection{Proof of Theorem 1 for $r\geq r_0$}

In this subsection we discuss the proof of Theorem 1 of the main paper for $r\geq r_0$. The proof follows essentially the same ideas of the proof in the case $r<r_0$ but there are some technical differences. Note that, the analysis of risk at the origin is unchanged because Lemma \ref{Lemma: Origin} holds for all $r$. So now, we analyze risk at non-origin points and basically prove equation \eqref{Eq: aux3}. As before, we use the decomposition of risk in equation \eqref{Eq: decomp}. Our strategy is the same, that is, showing that contribution of $\Exp D_{\theta \ld(\theta)}(Z)$ for one particular index $\ld(\theta)$ is dominant in $\Exp\log D_\theta(Z)$ and using a naive lower bound on $\Exp\log N_{\theta v}(Z)$ considering $\Exp N_{\theta \ln(\theta)}$ for one particular index $\ln(\theta)$. The choices of the indexes in this case are slightly different. Recall that each cluster $C_j$ of $\pi_{\sf C}[\eta,r]$ consists of only one point. The atoms are at $\mu_p=p\lamf$ for all $p\in\ZZ$. Without loss of generality, we only consider $\theta>0$. By Lemma \ref{Lemma: Negative}, which didn't depend on value of $r$, we can ignore all $D_{\theta p}$ with $p<0$. Suppose $\theta\in[\mu_l,\mu_{l+1})$ for $l\geq 0$. The contribution of $D_{\theta i}$ for all $i>l$ is negligible compared to $D_{\theta l}$ and the proof is exactly same as done in the beginning of Lemma \ref{Lemma: D}, c.f., equations \eqref{Eq: aux2}, \eqref{Eq: aux1}, \eqref{Eq: aux3} and \eqref{Eq: aux7}. The crucial difference from the sub-critical case arises now. We will see that, if $l\geq 1$, then unlike the sub-critical case $D_{\theta l}$ is not always the dominant term. Instead $D_{\theta,l-1}$ dominates $D_{\theta l}$ for some $\theta$ if $r>r_0$. To see this, note that,
\[
\Exp\log\Lp 1+\frac{D_{\theta,l-1}}{D_{\theta l}}\Rp=\Exp\log\Lp 1+c_{l-1}c_l^{-1}\exp\LP\lamf\Lp Z+\frac{\mu_l+\mu_{l-1}}{2}-\theta\Rp+\frac{\lame^2}{2}\Rp\RP.
\]
Hence, $\Exp D_{\theta l}(Z)$ dominates $\Exp D_{\theta, l-1}(Z)$ if $\lamf((\mu_l+\mu_{l-1})/2-\theta)+\lame^2/2\leq 0$, which simplifies to $\theta\geq\mu_l+\lamf/(2r)$. For $\theta\in[\mu_l,\mu_l+\lamf/(2r)]$, it can be shown that $D_{\theta, l-1}$ is dominant. Also note that for $l\geq 2$
\[
\Exp\log\Lp 1+\frac{D_{\theta,l-2}}{D_{\theta,l-1}}\Rp=\Exp\log\Lp 1+c_{l-2}c_{l-1}^{-1}\exp\LP\lamf\Lp Z+\frac{\mu_{l-2}+\mu_{l-1}}{2}-\theta\Rp+\frac{\lame^2}{2}\Rp\RP.
\]
Using $r\geq 1/2$
\[
\lamf\Lp\frac{\mu_{l-2}+\mu_{l-1}}{2}-\theta\Rp+\frac{\lame^2}{2}\leq\frac{\lame^2}{2}-\frac{3\lamf^2}{2}\leq 0.
\]
Hence, $D_{\theta, l-2}$ and the preceding $D_{\theta, i}$s are not dominant. Now if $D_{\theta l}$ is dominant, that is, $\theta\in[\mu_l+\lamf/(2r),\mu_{l+1})$ then using the naive lower bound $\Exp\log N_{\theta,v}(Z)\geq \Exp\log N_{\theta l}(Z)$ we get 
\[
\rho(\theta,\hat{p}_{\sf C}[\eta,r])=\frac{\theta^2}{2r}-\Exp\log N_{\theta,v}(Z)+\Exp\log D_\theta(Z)\leq\frac{\lamf^2}{2r}\Lp 1+o(1)\Rp.
\]
We skip the details of the proof because it's exactly similar to the case $\ld=\ln$ in Lemma~\ref{Lemma: DNCancel}. On the other hand, if $D_{\theta, l-1}$ is dominant, that is, $\theta\in[\mu_l,\mu_l+\lamf/(2r))$, then we use $\Exp\log N_{\theta l}(Z)$ as a lower bound of $\Exp\log N_{\theta,v}(Z)$. We end up with a quadratic in $\theta$ similar to equation \eqref{Eq: Quadratic}, which is nonnegative in $[\mu_l,\mu_l+2r\lamf]$. Since this interval covers the interval $[\mu_l,\mu_l+\lamf/(2r)]$ for $r\geq r_0$ we get the the above equation. 

We can also show that the cutoff $r_0=1/2$ is actually optimal. We discuss this briefly. Consider $r<1/2$ but suppose we still use $K=1$ so that the atoms are still at $\mu_p=p\lamf$ for $p\in\ZZ$. If $l\geq 0$, then the exact range of $\theta$ for which $D_{\theta l}$ is dominant is $[\mu_l+\lamf/(2r),\mu_l+\lamf(1+1/(2r))]$. If $r$ is small, this can be far from $\mu_l$, if $r$ is close to $1/2$ then it's close $[\mu_{l+1},\mu_{l+2}]$. On the other hand, we can show that the exact range where $N_{\theta p}$ is dominant is $[\mu_p-\lamf/2,\mu_p+\lamf/2]$. Now fix an $l\geq 0$. Consider all $\theta\in[\mu_l+\lamf/(2r),\mu_l+\lamf(1+1/(2r))]$. Since this interval has length $\lamf$, there exists $p>l$ such that, either $N_{\theta p}$ or $N_{\theta,p+1}$ is dominant for each $\theta$. Let $n=p-l$. Using these bounds in \eqref{Eq: decomp} to prove \eqref{Eq: nonorigin} we again get a quadratic as in \eqref{Eq: Quadratic}, which has roots
\begin{align*}
\alpha_n& =\mu_l+(1+r)n\lamf-\sqrt{n^2r^2+n^2r-nr-n+1},\mbox{ and}\\
\beta_n&=\mu_l+(1+r)n\lamf+\sqrt{n^2r^2+n^2r-nr-n+1}.
\end{align*}
If $n=1$, then $\alpha_1=\mu_{l+1}=\mu_l+\lamf$ and $\beta_1=\mu_l+(1+2r)\lamf$. If $n=2$, then $\alpha_2=\mu_l+\lamf(2(1+r)-\sqrt{(2r+1)^2-2(r+1)}$. For $r=1/2-\epsilon$ with very small $\epsilon>0$, we can check $\beta_1<\alpha_2<\mu_l+\lamf(1+1/2r)<\beta_2$. The small gap between $\beta_1$ and $\alpha_2$ makes the risk increase beyond threshold $\lamf^2/2r$. The cutoff $r_0=1/2$ is optimal because for $r=1/2$ we have the equality $\beta_1=\alpha_2=\mu_l+\lamf(1+1/(2r))$. Similarly, even for smaller values of $r$ one can show that such a gap always exists where the risk goes above the threshold $\lamf^2/(2r)$. These increments are of course order of $\Theta(\lamf^2)$ so that \eqref{Eq: nonorigin} fails to hold. This leads to the following result.

\begin{proposition} 
	If $\eta_n=s_n/n \to 0$, then for any $r < r_0=0.5$ and $\gamma \geq 1$,
	\[
	\limsup_{n \to \infty} \;\;\; \bigg\{\sup_{\theta \in \Theta_0[s_n]} \rho\big(\theta,\hat p_{\sf CL}[\eta_n,r; \gamma,  1]\big)\bigg\}\bigg/ R^*(\Theta_0[s_n]) \; >1.
	\]
	%$$\liminf_{n \to \infty}  B(\pi^{\sf C}_n)/ R^*(\Theta_0[s_n])=1.$$
\end{proposition}

\section{Proof of Theorem 2}
The Bayes risk of the multivariate cluster prior $B(\pi_n^{\sf C})=n B(\pi_{\sf C})$ and the univariate Bayes risk is given by $$ B(\pi_{\sf C})=(1-\eta_n) \rho(0,\hat p_{\sf C}[\eta_n,r]) + \frac{1-\eta_n}{2K}\sum_{i=1}^{\infty} \sum_{j=1}^K \eta_n^{|i|} \rho(\mu_{ij},\hat p_{\sf C}[\eta_n,r])$$
where $K$ is defined in eqn. (3) of the main paper. From the risk calculations in lemmas~2.1-2.3, it is clear that the first order asymptotic risk as $\eta_n \to 0$ can be reduced to just concentrating on the origin adjoining clusters $C_{\pm 1}$ and thereafter by symmetry: 
$$ B(\pi_{\sf C})= \frac{(1-\eta_n)\eta_n}{K} \sum_{j=1}^K  \rho(\mu_{1j},\hat p_{\sf C}[\eta_n,r]) (1+o(1))~.$$
Now, by lemma 2.2 and \eqref{Eq: decomp}, for each $1\leq j \leq K$ we have: 
$$  \rho(\mu_{1j},\hat p_{\sf C}[\eta_n,r]) =  \mu_{1j}^2/{(2r)} - \mathbb{E} \log  N_{\mu_{1j},v}(Z) +O(\lambda_{f,n})~ \text{ as } \eta_n \to 0~.$$
Also, following exactly the similar  asymptotic analysis as in lemma 2.2 abet now with $N_{\mu_{1j},v}(Z)$ we can  establish  $\mathbb{E} \log  N_{\mu_{1j},v}(Z)=(2v)^{-1}(\mu_{1j}^2  - \lambda^2_{f,n}) (1+o(1)).$ By construction, $\mu_{1j}  \geq  \lambda_{f,n}$ with strict equality only when $j=1$ and so each of the terms barring the first one has some positive contributions. Thus, $\rho(\mu_{11},\hat p_{\sf C}[\eta_n,r]) = \lambda_{f,n}^2/(2r)(1+o(1))$. For all $j > 1$, recalling  $\mu_{1j}/\lambda_{f,n}=(1+4r)^i \vee v^{-1/2} $ and $v=1/(1+r^{-1})$  we have,
$$\rho(\mu_{1j},\hat p_{\sf C}[\eta_n,r]) = 2^{-1} \lambda_{f,n}^2 \big\{1+r^{-1} - (1+4r)^{2i}\big\}_+ + O(\lambda_{f,n}) $$
where, the first term in the right side above is $0$ only when $j=K$. Thus, the maximal risk is only attained at $\mu_{11}=\lambda_{f,n}$. Thereafter, the risk decays and finally at $j=K$, the risk is negligible compared to the asymptotic minimax risk. Figure~1 shows the numerical evaluation of the risk of the cluster prior at the different support point of the first cluster. The figure shows the risk profile when $\eta_n= 10^{-15}$ which well captures the asymptotic analysis and the aforementioned decay in the risk function is evident from the figure. 

\begin{figure}[h]
\includegraphics[height=8cm, width=\textwidth]{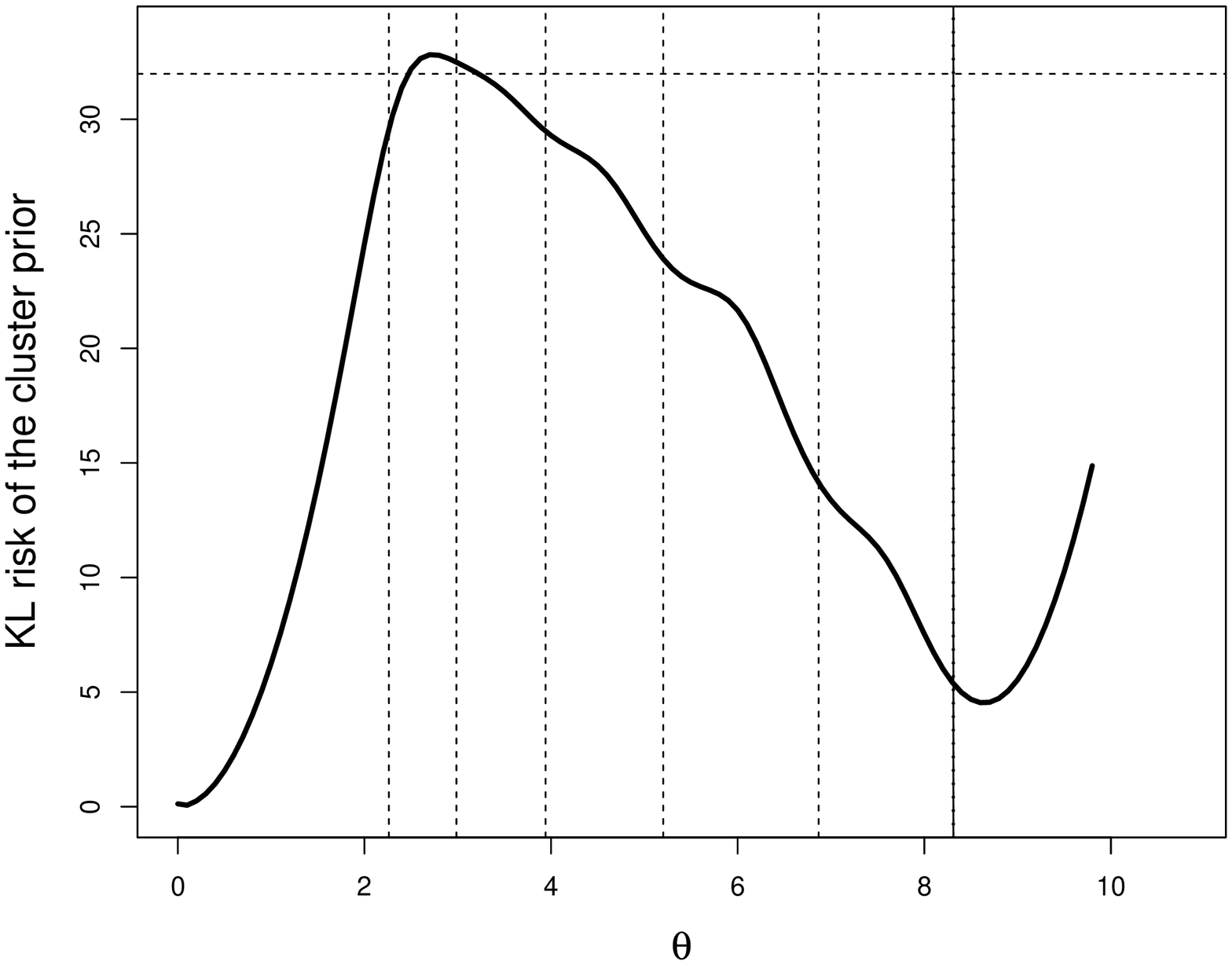}
\caption{Plot of the univariate predictive KL risk $\rho(\theta,\hat p_{\sf C}[\eta_n,r]$ as $\theta$ varies over the first cluster spanning $[0,\lambda_{e}(\eta_n)]$. Here, $\eta_n=10^{-15}$ and $r=0.08$. The horizontal line corresponds to the asymptotic theoretical limit $\lambda^2_f(\eta_n)/(2r)$. The dotted vertical lines denotes the location of the the  support points in cluster $C_1$ of $\pi_{\sf C}[\eta_n,r]$.}
\end{figure}

Noting that the multivariate minimax risk is $n \eta_n \lambda_{f,n}^2/(2r)(1+o(1))$ as $\eta_n \to 0$, the result follows from the above display. When $r  > r_0$, then $K=1$ and so, the above result directly imply $B(\pi_n^{\sf C})/R^*(\Theta_0[s_n]) \to 1$ as $n \to \infty$.  The condition $s_n \to \infty$ ensures that the prior concentrates on the parametric space $\Theta_0[s_n]$ defined in page 2 of the main paper (see Theorem 1B of \citet{Mukherjee-15} for details) and thus is least favorable in this case. 

\bibliography{pred-inf,append}